\documentclass[reqno,11pt]{amsproc}

\usepackage{amssymb}
\usepackage{amscd}
\usepackage{amsmath}

\newcommand{\ka}{{\kappa}} %
\theoremstyle{definition}
\newtheorem{dfn}{Definition}[section]

\newtheorem{remark}{Remark}

\newtheorem{thm}{Theorem}
\newtheorem{acknowledgement}{Acknowledgement}

\begin{document}
\keywords{Inverse problem, planar tree of strings,
Titchmarch-Weyl-matrix,
 dynamic Steklov-Poincar{\'e} operator,  leaf-peeling method}



\title[Two-velocity tree-like graph]
{On the inverse problem of the two-velocity tree-like graph}


\author[S. Avdonin]{Sergei Avdonin}%
\address{Department of Mathematics and Statistics, University of Alaska, Fairbanks,
AK 99775-6660, USA.}
\author[A. Choque Rivero]{Abdon Choque Rivero}
\address{Instituto de F\'{\i}sica y Matem\'aticas, Universidad Michoacana de San
  Nicol\'as de Hidalgo, C.P. 58040 Morelia, Mich., M\'exico.} 
\author[G. Leugering]{Guenter Leugering}
\address{Department Mathematik, Lehrstuhl Angewandte Mathematik II, Cauerstr.
 11 91058 Erlangen Germany.}
\author[V. Mikhaylov]{Victor Mikhaylov}
\address{St. Petersburg Department of V.A. Steklov Institute of Mathematics
of the Russian Academy of Sciences, 27 Fontanka, St-Petersburg,
Russia { and St. Petersburg State University, Faculty of
Physics, Ulyanovskaya str. 3, 198504 St-Petersburg,
Russia. }}

\begin{abstract}
In this article the authors continue the discussion in \cite{ALM}
 about inverse problems for second order elliptic and hyperbolic equations
  on metric trees from boundary measurements. In the present paper we prove
  the identifiability of varying densities of a planar tree-like network of
   strings along with the complete information on the graph, i.e. the lengths
    of the edges, the edge degrees and the angles between neighbouring edges.
    The results are achieved using the Titchmarch-Weyl function for the spectral
     problem and the Steklov-Poincar{\'e} operator for the dynamic wave equation
      on the tree. The general result is obtained by a peeling argument which reduces
       the inverse problem layer-by-layer from the leaves to the clamped
       root of the tree.
\end{abstract}
\maketitle                   





\section{Introduction}
Partial differential equations on networks become 
increasingly important in applied sciences and engineering. Starting
on the atomistic scale, the so-called quantum graphs play a
prominent role for physicists and chemists. On the nanoscale,
chemical engineers are interested in mechanical properties of
nanorods. On the mesoscale  fiber networks modeled by strings and
beams are crucial for tissue engineering, whereas on the macroscale
grid structures involving pipes, plates and shells appear in
abundance in civil and mechanical engineering. Clearly, the modeling
of networks of Schr\"odinger equations \cite{kuch,ks1} or those for
elastic strings, beams, plates and shells
\cite{LLS1994,LagneseLeugering2004} is by now a classical subject
and is well established. This applies also to problems of optimal
control and optimization of quantum graphs \cite{A}, string and beam
networks, mostly, however, in the linear range. See e.g.
\cite{AvdoninIvanov1995,LLS1994,LagneseLeugeringSchmidt1994,
LagneseLeugeringSchmidt1993,LagneseLeugering2004,DagerZuazua2006,AM},
where problems of exact controllability and uniform stabilizability
have been considered for tree-like networks.
 As for the modeling and
controllability of tree-like networks of nonlinear strings see
\cite{LeugeringSchmidt2011}. Similar work for nonlinear Timoshenko
beams and Cosserat rods is under progress.

 It is important to note
that most controllability and uniform stabilizability results known
in the literature rely on a single and simple Dirichlet node and
controls on all remaining leaves of the tree; see exceptions
in\cite{LeugeringZuazua1999,DagerZuazua2006}, where a tree is
controlled from the root, though, under severe number theoretic
restrictions on the size of the strings
 involved.
 In
\cite{LeugeringSokolowski2009} planar networks of elastic strings
have been investigated under topology changes. There a sensitivity
analysis with respect to the nodal configuration was developed.
Together with thickness optimization for beam networks, this leads
to design problems for such networks.

 For inverse problems which are
strongly related to design problems governed by partial or even
ordinary differential equations on metric graphs, the literature is
far less developed, in particular for planar  problems. On the other
hand, such problems for quantum graphs are widely studied; see
\cite{AK,ALM,AMR,B,BV,BW,FY,KuSt,Y} for a not exhaustive account on
the literature. To the best knowledge of the authors, \cite{ALM} is
the first article in which the inverse problem for planar elastic
string tree networks has been investigated mathematically. The model
for the linear elastic strings relied on constant coefficients for
the equations on each edge. It was shown that these coefficients and
the full topological information of the graph, namely the edge
degrees, the lengths of the edges and the angles between the edges
incident at a given multiple node could be recovered from the
(reduced) dynamical response operator or Titchmarsh-Weyl matrix. In
this article the authors extend those results to tree-like networks
of planar elastic strings with spatially varying coefficients. This
step is very important in failure detection of such networks.

More precisely, we consider the in-plane motion of elastic strings
forming  to a tree-like network. Mathematically this means that a
 two-velocity wave equation for a two-component vector displacement is assumed to hold on each edge of the
tree, while transmission conditions at the vertices of the graph
couple the edgewise-defined differential equations in order to form
a differential system on the entire metric graph.
 We
investigate the inverse problem of recovering  the physical
properties, i.e.
  the (spatially varying) densities and lengths of each string, and also
 the topology of the tree and
 the angles between branching
edges. It is shown that the inverse problem can be uniquely solved
by applying measurements at all, or at all but one boundary
vertices.

 Let $\Omega$ be a finite connected compact planar graph without cycles, i.e. a tree,
 with edges $E=\{e_1, \dots, e_N\} = \{e_i | i \in \mathcal{I}\}$, where
 $\mathcal{I}:=\{1,\dots,N\}$, connected at the vertices (nodes) $V =
 \{v_1,\dots,v_{N+1}\}$. Every edge $e_i \in E$ is identified with an interval
$[\beta_{2i-1},\beta_{2i}]$ of the real line and is equipped with
two  positive functions (densities): $\varrho_i,\,\mu_i\in C^2(e_i).$

The displacements at each point of the graph are denoted by $r^i(x)
\in \textbf{R}^2$. As each edge carries along an individual
coordinate system $\textbf{e}_i, \textbf{e}_i^{\perp}$, the
displacement decomposes as follows:
\begin{equation}
r^i(x) := u^i(x)\textbf{e}_i + w^i(x)\textbf{e}_i^{\perp},
\end{equation}
where $u^i(x)$ is the longitudinal (tangential) displacement and
$w^i(x)$ the vertical (in-plane) or normal displacement at the
material point $x \in [\beta_{2i-1},\beta_{2i}]$.

    Given a node $v_J$ we define
    $\mathcal{I}_J:= \{i \in \mathcal{I} | e_i\mbox{ is incident at
    }v_J\}$,
    the incidence set, and $d_J = |\mathcal{I}_J|$ the edge degree of $v_J$.
     The set $\mathcal{J} = \{J | v_J \in V\}$
    of node indices splits into $\mathcal{J}_S$ and $\mathcal{J}_M$ which
    correspond to simple and multiple nodes according
    to $d_J=1$ and $d_J > 1$, respectively. The set of simple nodes
    constitutes the boundary of the graph and is denoted by $\Gamma.$

    In the global Cartesian coordinate system, one can represent each edge by a rotation matrix
    \begin{equation*}
        S_\alpha =
        \left(
            \begin{matrix}
                \cos \alpha & -\sin \alpha \\
                \sin \alpha & \cos \alpha
            \end{matrix}
        \right),
    \end{equation*}
    where $\textbf{e} = (\cos \alpha, \sin \alpha)^T$. In fact,
    it is the global coordinate system that
     we will use, as we are going to identify the angles $\alpha_{ij}$
      between two branching edges.

\section{Spectral and  dynamical settings}

The space of real vector valued square integrable functions on the
graph $\Omega$ is denoted by $L_2(\Omega) =
\bigoplus_{i=1}^{N}L_2(e_i, \mathbb{R}^2).$ For the element $U \in
L_2(\Omega)$ we write
$$
U = \{u,w\} =
        \left\{
            \left(
                \begin{matrix}
                    u^i \\
                    w^i
                \end{matrix}
            \right)
\right\}_{i=1}^N\ ,\ u^i, w^i \in L_2(e_i).
$$
    Denote by $\alpha_{ij}$ the angle between
     two edges $e_i$ and $e_j$ counting from $e_i$
      counterclockwise, and introduce the matrices
    \begin{equation}
        D_i (x)=
        \left(
            \begin{matrix}
                \frac{1}{\varrho_{i}(x)} & 0 \\
                0 & \frac{1}{\mu_{i}(x)}
            \end{matrix}
        \right), \ i \in \mathcal{I}. \label{2.4}
    \end{equation}
The continuity condition at the internal vertex $v_J$ reads
\begin{equation}
\left(\begin{matrix}
u^i(v_J) \\
w^i(v_J)
\end{matrix}
\right) = S_{ij} \left(
\begin{matrix}
u^j(v_J) \\
w^j(v_J)
\end{matrix}
\right), \ \ i,j \in \mathcal{I}_J, \ \ S_{ij}:=S_{\alpha_{ij}}.  \label{2.5}
\end{equation}
The second condition (balance of forces) at the internal vertex
$v_J$ is
\begin{equation}
\sum_{j \in \mathcal{I}_J} S_{ij} D_{j}(v_J) \left(
\begin{matrix}
u_x^j(v_J) \\
w_x^j(v_J)
\end{matrix}\right) =0. \label{2.6}
\end{equation}
It is easy to check that if  (\ref{2.6}) is satisfied for some
 $i \in \mathcal{I}_J,$ then it is valid for any  $i \in \mathcal{I}_J.$
The compatibility conditions (\ref{2.5}), (\ref{2.6})  are called the
 Kirchhoff--Neumann (KN) matching conditions.

We put $\{\phi, \psi\} = \left\{ \left( \begin{matrix} \phi^i \\
\psi^i \end{matrix} \right) \right\}_{i=1}^N \in L_2(\Omega)
    $,\, $\phi^i,\, \psi^i\in H^2(e_i)$
     and associate the following spectral problem to the graph.
\begin{dfn}
By $\textbf{S}$ we denote the spectral problem on the graph. It is
described by the equation on the edges
\begin{align}
\left.\begin{array}{cc}
-\phi_{xx}^i & = \lambda \varrho_{i}(x)\phi^i \\
-\psi_{xx}^i & = \lambda \mu_{i}(x)\psi^i
\end{array}\right\}
\ \ x \in e_i, \ \ i \in \mathcal{I}; \label{2.7}
\end{align}
\begin{equation}
        \{\phi, \psi\} \ \  \mbox{satisfies the KN conditions}\ \ (\ref{2.5}), (\ref{2.6}) \mbox{ at all internal vertices} \ \
        v_J, J \in \mathcal{J}_M, \label{2.8}
\end{equation}

        \begin{equation}\{\phi,\psi\} = 0\ \  \mbox{on the boundary} \ \
        \Gamma. \label{2.9}\end{equation}
 The last equality means that $\phi^i(v_J) = \psi^i(v_J)=0$, for $i \in \mathcal{I}_J$, $J \in \mathcal{J}_S$.
\end{dfn}

Along with the spectral description, we consider the dynamical
system, described by the two-velocity problem on the each edge of
the graph:
\begin{eqnarray}
\label{2.10} \left.\begin{array}{cc} \varrho_{i}(x)\,u_{tt}^i -
u_{xx}^i = 0\\
\mu_{i}(x)\, w_{tt}^i - w_{xx}^i =0 \end{array}\right\} \quad  t >
0,\ \ x \in e_i,\ \  i \in \mathcal{I}.
\end{eqnarray}
We assume that $|\mathcal{J}_S|=m$. By $\mathcal{F}_\Gamma^T =
L_2([0,T], \mathbb{R}^{2m})$ we denote the space of controls acting
on the boundary of the tree. For the element $F \in
\mathcal{F}_\Gamma^T$ we write $F = \{f,g\} = \left\{ \left(
\begin{array}{c}
f^i \\
g^i
\end{array}
\right) \right\}_{i=1}^m\ ,\ f^i, g^i \in L_2(0,T).$

We deal with Dirichlet boundary conditions
\begin{equation}
\{u,w\} = \{f,g\},\ \ \mbox{ on }\Gamma \times [0,T], \label{2.12}
\end{equation}
where $f,g \in \mathcal{F}_\Gamma^T,$ and the last equality means that $u^i(v_J) =f^i$, $w^i(v_J)
= g^i$, for $i \in \mathcal{I}_J$, $J \in \mathcal{J}_S$.
\begin{dfn}
By \textbf{D} we denote the dynamical problem on the graph $\Omega$,
described by the Eqs.  \eqref{2.10} which satisfies the
compatibility conditions \eqref{2.5}, \eqref{2.6} at all internal
vertices for all $t > 0$, Dirichlet boundary condition \eqref{2.12}
and zero initial conditions $\{u(\cdot, 0) , w(\cdot , 0)\} =
\{0,0\}$, $\{u_t(\cdot , 0) , w_t(\cdot , 0)\} = \{0,0\}$.
\end{dfn}
It is known that for any $T>0$, there exists a unique generalized
solution of this problem such that $\{u , w\} \in C([0,T] ;
L^2(\Omega))$ if $F \in \mathcal{F}_\Gamma^T$; see, e.g.
\cite{LLS1994}.

\section{Inverse  spectral  and  dynamical  problems}
We use the Titchmarsh-Weyl (TW) matrix function as the data for the
inverse spectral problem. Let $\lambda \notin \mathbb{R}$. We
introduce  two solutions of \eqref{2.7}, \eqref{2.8}, $\{\phi^1,
\psi^1\}$ and  $\{\phi^2, \psi^2\}$ with the following boundary
conditions
    \begin{equation}
        \{ \phi^1 , \psi^1 \} =
        \left(
            \begin{matrix}
                (0,0) \\
                \cdots \\
                (1,0) \\
                \cdots \\
                (0,0)
            \end{matrix}
        \right), \ \
        \{ \phi^2 , \psi^2 \} =
        \left(
            \begin{matrix}
                (0,0) \\
                \cdots \\
                (0,1) \\
                \cdots \\
                (0,0)
            \end{matrix}
        \right)\ \ \mbox{ on }\Gamma,\label{3.1}
    \end{equation}
where nonzero elements are located at the $i$-th row. Then the TW
matrix ${\bf M}(\lambda$) is defined
     as ${\bf M}(\lambda) = \{M_{ij}(\lambda)\}_{i,j=1}^m$ where each $M_{ij}(\lambda)$ is a $2 \times 2$ matrix defined by
    \begin{equation}
        M_{ij}(\lambda) =
        \left(
            \begin{matrix}
                \phi_x^1(v_j,\lambda) & \psi_x^1(v_j,\lambda) \\
                \phi_x^2(v_j,\lambda) & \psi_x^2(v_j,\lambda)
            \end{matrix}
        \right), \ \ 1 \le i,j \le m. \label{3.2}
    \end{equation}
     Let us consider problem \eqref{2.7} with the nonhomogeneous Dirichlet boundary condition
    \begin{equation}
        \{ \phi, \psi \} = \{ \zeta, \nu \} \mbox{ on }\Gamma, \label{3.3}
    \end{equation}
where $\{ \zeta, \nu \} \in \mathbb{R}^{2m}$, and let $\{\phi,
\psi\}$ be the solution to
    \eqref{2.7}, \eqref{2.8}, \eqref{3.3}.
     The Titchmarsh-Weyl matrix connects the values of the solution $\{ \phi, \psi \}$ on the boundary and
      the values of its derivative $\{ \phi_x , \psi_x \}$  on the
      boundary:
    \begin{equation}
        \{\phi_x , \psi_x \} = {\bf M}(\lambda)\{ \zeta , \nu \}\ \mbox{ on }
        \Gamma. \label{3.4}
    \end{equation}
We set up the inverse spectral problem as follows: Given the TW
matrix function ${\bf M}(\lambda)$,
     $\lambda \notin \mathbb{R}$, we have to recover the graph
     (lengths of edges, connectivity and angles between edges)
     and parameters of the system \eqref{2.7}, i.e. the set of the
     densities $\{\varrho_{i},
     \mu_{i}\}_{i=1}^N$.

 Let $\{u,w\}$ be the solution to the problem \textbf{D} with the boundary control $\{f,g\} \in \mathcal{F}_\Gamma^T$.
 We introduce the \textit{dynamical response operator} (the dynamical Dirichlet-to-Neumann map)
 to the problem \textbf{D} by
    \begin{equation}
        R^T \{f,g\}(t) = \left.\{u_x(\cdot, t) , w_x(\cdot, t)\}\right|_\Gamma , \ \   t \in [0,T].
         \label{3.5}
    \end{equation}
    The response operator has the form of a convolution:
    \begin{equation}
        \left( R^T \{f,g\} \right) (t) = (\textbf{R} \ast \{f,g\})(t)=\int_0^t \textbf{R}(t-s)\{f,g\}(s)\,ds, \ \ t \in
        [0,T], \label{3.6}
    \end{equation}
    where $\textbf{R}(t) = \{R_{ij}(t)\}_{i,j=1}^m$ is the response matrix
     function, and each entry $R_{ij}(t)$ is a $2 \times 2$ matrix.
    The entries $R_{ij}(t)$ are defined by the following procedure. We set up two
    dynamical problems defined by Eqs.
     \eqref{2.10}, \eqref{2.5}, \eqref{2.6}
     and the boundary conditions given by
    \begin{equation*}
        \{U^1(\cdot , t) , W^1(\cdot , t)\} =
        \left(
            \begin{matrix}
                (0,0) \\
                \cdots \\
                (\delta(t),0) \\
                \cdots \\
                (0,0)
            \end{matrix}
        \right),
        \{U^2(\cdot , t) , W^2(\cdot , t)\} =
        \left(
            \begin{matrix}
                (0,0) \\
                \cdots \\
                (0,\delta(t)) \\
                \cdots \\
                (0,0)
            \end{matrix}
        \right) \mbox{ on } \Gamma,
    \end{equation*}
 where $\delta(t)$ is the delta distribution.

    In the above formulas, the only nonzero row is the $i$-th. Then
    \begin{equation*}
        R_{ij}(t) =
        \left(
            \begin{matrix}
                U_x^1(v_j, t) & W_x^1(v_j, t) \\
                U_x^2(v_j, t) & W_x^2(v_j, t)
            \end{matrix}
        \right).
    \end{equation*}
    Therefore, in order to construct the entries of $\textbf{R}$, we need to set up the boundary condition at the $i$-th
    boundary node in the first and second channels, while having other boundary nodes fixed
    (impose homogeneous Dirichlet conditions there) and measure
     the response at $j$-th boundary node in the first and second channels.

We set up the dynamical inverse problem as follows: given the
response operator
     $R^T$ \eqref{3.5} (or what is equivalent, the matrix $\textbf{R}(t),$  $t \in [0,T]$)
      for large enough $T$, we recover the graph (lengths of edges, connectivity and angles
       between edges)
       and the parameters
    of the dynamical system \eqref{2.10}, i.e. speeds of the wave propagation on the
     edges. The connection between the spectral and the dynamical data is known (see, for example,  \cite{AK,ALM,AMm,AMR,KKLM})
      and was used for inverse spectral and dynamical problems.
    Let $\{f,g\ \} \in \mathcal{F}_\Gamma^T \cap (C_0^{\infty}(0, +\infty))^{2m}$ and
    \begin{equation*}
        \widehat{\{f,g\}}(k) := \int_{0}^{\infty} \{ f(t) , g(t) \} e^{ikt}dt
    \end{equation*}
    be its Fourier transform. Eqs. \eqref{2.10}, and \eqref{2.7}
    are clearly connected: Going formally in \eqref{2.10}
     over to the Fourier transform, we obtain \eqref{2.7} with $\lambda = k^2$.
     It is not difficult to check  that the response matrix function ${\bf R}(t)$
      and
     TW matrix function $\textbf{M}(\lambda)$ (Nevanlinna type matrix
     function) are connected by the same transform:
    \begin{equation}
    \label{Sp_Dyn}
    \textbf{M}(k^2) = \int_{0}^{\infty} {\textbf R}(t) e^{ikt}dt
    \end{equation}
    where this equality is understood in a weak sense.



\section{Inverse problem on a star graph}
 In this section we consider a graph consisting of $n$ edges $e_i=[0,l_i]$,
$i=1,\ldots,n$, connected at the internal vertex identified with the
point $x=0;$ the boundary vertices are identified with  the
endpoints $l_i.$ We will show that to recover the graph and the
parameters of the system it is sufficient to use the response
operator or the TW matrix, associated with all or all but one
boundary vertices; we denote this vertex by $v_n$.

On each edge we have the following equations
\begin{equation}
\label{Star_eqn}
\left.\begin{array}{cc}\varrho_{i}(x)\,u_{tt}^i- u_{xx}^i = 0\\
\mu_{i}(x)\, w_{tt}^i - w_{xx}^i = 0
\end{array}
\right\} \quad  t > 0, \, x \in (0,l_i),\ \ i\in\{1,,\ldots,n\}.
\end{equation}
We introduce the notations
\begin{equation*}
D_i=\left(\begin{matrix}
\frac{1}{\varrho_{i}(0)} & 0 \\
0 & \frac{1}{\mu_{i}(0)}
\end{matrix} \right), \quad
B_i =\left(\begin{matrix}
\frac{1}{\sqrt{\varrho_{i}(0)}} & 0 \\
0 & \frac{1}{\sqrt{\mu_{i}(0)}}
\end{matrix}\right),
\quad i =1,\ldots,n,
\end{equation*}
and write  the KN conditions at the internal vertex in the form
\begin{equation}
\left(\begin{matrix}
u^i(0,t) \\
w^i(0,t)
\end{matrix}\right) = S_{ij}
\left(\begin{matrix}
u^j(0,t) \\
w^j(0,t)\end{matrix} \right),\quad
i,j\in\{1,\ldots,n\},\label{Star_cont}
\end{equation}
\begin{equation}
\sum_{j=1}^n S_{ij} D_j \left( \begin{matrix}
u_x^j(0,t) \\
w_x^j(0,t)
\end{matrix}\right) =0. \label{Star_Kirhg}
\end{equation}

\noindent {\bf Problem:} Given ${\bf M(\lambda)},$   an $(n-1)\times
(n-1)$ matrix whose entries are $2\times 2$ matrices $M_{ij}(\lambda), \; i,j=1,\ldots,n-1,$ find
$\varrho_i(x)$, $\mu_i(x)$, $l_i$, $i=1,\ldots,n$ and
$\alpha_{i,j}$, $i,j=1,\ldots,n$.

\noindent {\bf Solving the inverse problem on a star graph:}

{\bf Step 1.} We fix $i, \, i\not= n,$ and set up the initial
boundary value problem: On the edges and at the internal vertex the
conditions (\ref{Star_eqn})--(\ref{Star_Kirhg}) hold, and the
initial conditions are zero and the boundary conditions  are
\begin{equation*}
u^i(l_i,t) = \delta(t), \  w^i(l_i,t) = 0, \ u^j(l_j,t) =w^j(l_j,t)
= 0,\quad j\not=i.
\end{equation*}
 We introduce the notations
\begin{eqnarray}
\tau_i(x) = \int_{0}^{x}\sqrt{\varrho_i(s)}ds,
\ \ \nu_i(x) = \int_{0}^{x}\sqrt{\mu_i(s)}ds, \label{notations}\\
L_i = \int_{0}^{l_i} \sqrt{\varrho_i(x)}dx, \ M_i
=\int_{0}^{l_i}\sqrt{\mu_i(x)}dx, \ \, i=1, \ldots, n.\notag
\end{eqnarray}
For  $t < L_i $, the solution of the initial boundary value problem
has the form $$
                u^i(x,t) = \left[\frac{\varrho_1(x)}{\varrho_1(l_1)}\right]^{-\frac{1}{4}}
                \delta \left(t-\int_x^{l_1} \sqrt{\varrho_1(s)}ds \right)+h_0(x,t), \ \   w^i = u^j = w^j = 0, \ j \neq i. $$
        For $
        L_i < t < L_i + \epsilon, \;$ $ \epsilon < \min_{i=1,\ldots,n} \{L_i,M_i\},$
the solution of the initial boundary value problem has the form
\begin{eqnarray} \label{sta}
&u^i(x,t)
=\left[\frac{\varrho_i(x)}{\varrho_i(l_i)}\right]^{-\frac{1}{4}}
\delta \left(t-\int_x^{l_i} \sqrt{\varrho_i(s)}ds \right)
+ a^i_i\left[\frac{\varrho_i(x)}{\varrho_i(0)}\right]^{-\frac{1}{4}} \delta(t - \tau_i(x) -  L_i) + h_i(x,t), \\
&w^i(x,t) = b^i_i \left[\frac{\mu_i(x)}{\mu_i(0)}\right]^{-\frac{1}{4}} \delta(t-\nu_i(x)-L_i) + q_i(x,t), \\
&u^j(x,t) = a^i_j \left[\frac{\varrho_j(x)}{\varrho_j(0)}\right]^{-\frac{1}{4}} \delta(t-\tau_j(x)  -L_j) + h_j(x,t),\quad j\not=i, \\
&w^j(x,t) = b^i_j
\left[\frac{\mu_j(x)}{\mu_j(0)}\right]^{-\frac{1}{4}}
\delta(t-\nu_j(x) -L_j) + q_j(x,t),\quad j\not=i.  \label{sta1}
\end{eqnarray}
Here $h_j(x,t)$ and $q_j(x,t)$ are continuous functions, and $a^i_j,
\, b^i_j$, $j=1,\ldots,n$, are real numbers. We introduce the
notations
$$
\ka_i:=
\left[\frac{\varrho_i(0)}{\varrho_i(l_i)}\right]^{-\frac{1}{4}}, \ \
\gamma_i:= \left[\frac{\mu_i(0)}{\mu_i(l_i)}\right]^{-\frac{1}{4}}
$$
and make use of the continuity condition (\ref{Star_cont})  to get (equating the coefficients of $\delta$)
\begin{equation}
\label{cont1}
\left(\begin{matrix} \ka_i + a^i_i \\
b^i_i \end{matrix}\right) = S_{ij} \left(
\begin{matrix} a^i_j \\
b^i_j \end{matrix}\right).
\end{equation}
From the force balance condition (\ref{Star_Kirhg}) we obtain (equating the coefficients of $\delta'$)
\begin{equation}
\label{kir1}
B_i \left(\begin{matrix}-\ka_i + a^i_i \\
b^i_i
\end{matrix}\right) +
\sum_{j=1,\,j\not= i}^n S_{ij} B_j\left(\begin{matrix} a^i_j \\
b^i_j
\end{matrix}\right) = 0, \quad j\not=i.
\end{equation}

{\bf Step 2.}  Similarly we consider the initial boundary value
problem (\ref{Star_eqn})--(\ref{Star_Kirhg}) with zero initial
condition and the following boundary conditions:
\begin{equation*}
u^i(l_i,t) = 0, \  w^i(l_i,t) = \delta(t), \ u^j(l_j,t) =w^j(l_j,t)
= 0,\quad j\not=i.
\end{equation*}
The solution for $M_i < t < M_i + \epsilon$ has the form
\begin{eqnarray}  \label{star}
&u^i(x,t) =\widetilde a^i_i
\left[\frac{\varrho_i(x)}{\varrho_i(0)}\right]^{-\frac{1}{4}}
\delta(t-\tau_i(x)-M_i)+\widetilde h_i(x,t),\\
&w^i(x,t) = \left[\frac{\mu_i(x)}{\mu_i(l_i)}\right]^{-\frac{1}{4}}
\delta \left(t-\int_x^{l_i} \sqrt{\mu_i(s)}ds \right) +\widetilde
b^i_i
\left[\frac{\mu_i(x)}{\mu_i(0)}\right]^{-\frac{1}{4}}\delta(t-\nu_i(x)-M_i)+
\widetilde q_i(x,t), \\
&u^j(x,t) = \widetilde a^i_j
\left[\frac{\varrho_j(x)}{\varrho_j(0)}\right]^{-\frac{1}{4}}
\delta(t-\tau_j(x)-M_j) + \widetilde h_j(x,t),\quad j\not=i, \\
&w^j(x,t) = \widetilde b^i_j
\left[\frac{\mu_j(x)}{\mu_j(0)}\right]^{-\frac{1}{4}}
\delta(t-\nu_j(x) -M_j) + \widetilde q_j(x,t),\quad j\not=i .  \label{star1}
\end{eqnarray}
Here $\widetilde h_j(x,t)$ and $\widetilde q_j(x,t)$ are continuous
functions, and  $\widetilde a^i_j, \, \widetilde b^i_j$,
$j=1,\ldots,n$, are real numbers.

From the matching conditions (\ref{Star_cont}),
(\ref{Star_Kirhg}) at the internal vertex, we obtain
\begin{equation}
\label{cont2}
\left(\begin{matrix} \widetilde b^i_i \\
\ka_i + \widetilde a^i_i \end{matrix}\right) = S_{ij} \left(
\begin{matrix} \widetilde b^i_j \\
\widetilde a^i_j \end{matrix}\right), \ \
B_i \left(\begin{matrix} \widetilde b^i_i \\
-\gamma_i + \widetilde a^i_i
\end{matrix}\right) +
\sum_{j\not=i} S_{ij} B_j\left(\begin{matrix} \widetilde b^i_j \\
\widetilde a^i_j
\end{matrix}\right) = 0, \quad j\not=i.
\end{equation}

 {\bf Step 3.}  At this point we use the boundary control (BC) method
 which is based on deep connections between controllability and identifiability of dynamical systems.
  This method and its connections with other approaches to inverse problems are described in detail in \cite{AMm,BIP07,B1,BM}.

From  ${\bf M(\lambda)}$ we find ${\bf R(t)}$;  see equation
\eqref{Sp_Dyn}.
 The characteristic feature of the BC method is its {\sl locality}: we can recover parameters of the edges as of independent intervals.
Eqs. \eqref{sta}--\eqref{sta1}, \eqref{star}--\eqref{star1}
demonstrate that we can apply the BC method for one interval (see,
e.g. \cite{AMm,B1,BM}) to find $\varrho_{i}(x)$ from the response
matrix functions $\{R_{ii}(t)\}, \; 0<t<2L_i,$ and to find
$\mu_{i}(x)$ from $\{R_{ii}(t)\}, \; 0<t<2M_i.$ Furthermore,
 from  $\left\{R_{ii}(t)\right\}$ we can recover the
transmission and reflection coefficients $a^i_i,$ $b^i_i$,
$\widetilde a^i_i,$ $\widetilde b^i_i$.

{\bf Step 4.} In this step we recover $B_n$ and $\alpha_{ij}.$

\noindent We introduce the notations $\xi_i=\ka_i+a^i_i,$
$\eta_i=b^i_i$, $\widetilde \xi_i=\widetilde b^i_i,$
$\widetilde\eta_i=\gamma_i+\widetilde a^i_i$. Then from
(\ref{cont1}), (\ref{kir1}) we deduce
\begin{equation}
\label{kir1_1}
\sum_{j=1,\,j\not= i}^n S_{ij} B_j \left(S_{ij}\right)^{-1}\left(\begin{matrix} \xi_i \\
\eta_i
\end{matrix}\right) +B_i \left(\begin{matrix}\xi_i \\
\eta_i
\end{matrix}\right)= 2B_i\left(\begin{matrix}\ka_i \\
0
\end{matrix}\right), \quad j\not=i.
\end{equation}
The Eqs. (\ref{cont2}) yield
\begin{equation}
\label{kir2_1}
\sum_{j=1,\,j\not= i}^n S_{ij} B_j \left(S_{ij}\right)^{-1}\left(\begin{matrix} \widetilde\xi_i \\
\widetilde\eta_i
\end{matrix}\right) +B_i \left(\begin{matrix}\widetilde\xi_i \\
\widetilde\eta_i
\end{matrix}\right)= 2B_i\left(\begin{matrix} 0 \\
\gamma_i
\end{matrix}\right), \quad j\not=i.
\end{equation}
Since the matrices in the LHS  of (\ref{kir1_1}) and  (\ref{kir2_1}) are positive,
we necessarily have $\left(\begin{matrix}\xi_i \\
\eta_i
\end{matrix}\right)\not= \left(\begin{matrix}\widetilde\xi_i \\
\widetilde \eta_i
\end{matrix}\right)$.
Thus  (\ref{kir1_1}), (\ref{kir2_1}) completely determine the
matrixes $A_i$:
\begin{equation}
\label{d3} A_i=\sum_{j=1}^n S_{ij} B_j (S_{ij})^{-1},
\quad i=1\ldots n-1.
\end{equation}
In (\ref{d3}) we do not know the $n-1$ angles between the edges and
the matrix $B_n$.

Multiplying (\ref{d3}) by $S_{ki}$ from the left and by
$(S_{ki})^{-1}$ from the right and using that
$S_{ki}=(S_{ik})^{-1}$, $S_{ki}S_{ij}=S_{kj}$, we obtain
\begin{equation} \label{S1}
A_k=S_{ki}A_i(S_{ki})^{-1} \ \ {\rm{or}} \ \  A_k S_{ki} =S_{ki}A_i.
\end{equation}
The angle $\alpha_{ki}$ can now be found from this relation.
Repeating this procedure for various $k$, $i$ we can determine all
angles $\alpha_{ki},\,$  $k,i\not=n$. Then we can use any of the
Eqs. (\ref{d3}) to determine $B_n$ and $\alpha_{in}$. Indeed,
 (\ref{d3}) can be written as
\begin{equation} \label{S2}
S_{in} B_n (S_{in})^{-1}=C_i  \ \ {\rm{or}} \ \ B_n= (S_{in})^{-1}C_i S_{in}
\end{equation}
with some known matrix $C_i$. Using the invariants
\begin{equation*}
\frac{1}{\sqrt{\rho_n(0)}}+\frac{1}{\sqrt{\mu_n(0)}}=\operatorname{tr}{C_i},\quad
\frac{1}{\sqrt{\rho_n(0)}}\frac{1}{\sqrt{\mu_n(0)}}=\det{C_i},
\end{equation*}
we determine the matrix $B_n$ and then the angle $\alpha_{in}$.

{\bf Step 5.} The leaf-peeling method.

In this step we reduce our inverse problem to the problem on a
 new tree, $\widetilde\Omega,$ consisting of  one edge,
$e_n=(0,l_n),$ and boundary vertices identified with the endpoints
of this interval.
 Our goal is to obtain the matrix element $\widetilde
M_{00}$ of the TW matrix for $\widetilde\Omega$ associated with the
``new" boundary, the vertex identified with $x=0.$

 We denote by $\{\Phi,\Psi\}=\left\{ \left(
\begin{matrix} \phi^i \\ \psi^i \end{matrix} \right)
\right\}_{i=1}^n$  the solution to the following problem on
the graph  $\Omega$:
\begin{align}
\left. \begin{array}{cc}
-\phi_{xx}^i & = \lambda \varrho_{i}(x)\phi^i \\
- \psi_{xx}^i & = \lambda \mu_{i}(x)\psi^i
\end{array}\right\}
\ \ x \in e_i, \ \ i=1,\ldots,n, \label{Star_sp_eqn}
\end{align}
with the KN conditions at the internal vertex and  the boundary
conditions
\begin{equation}
\label{Star_BCN} \psi^1(l_1,\lambda)=\zeta,\,\,
\phi^1(l_1,\lambda)=\nu,\,\,\psi^i(l_i,\lambda)=\phi^i(l_i,\lambda)=0,\quad
2\leqslant i \leqslant n.
\end{equation}
We note that  $\{\Phi,\Psi\}$
 solves the following Cauchy problem  on the edge $e_1:$
\begin{eqnarray*}
-\phi_{xx}^1=\lambda \varrho_1(x)\phi^1,\quad
-\psi_{xx}^1=\lambda \mu_1(x)\psi^1,\quad x\in (0,l_1), \\
\begin{pmatrix}
\phi^1(l_1) \\
\psi^1(l_1)
\end{pmatrix}=
\begin{pmatrix}
\zeta \\
\nu
\end{pmatrix}, \ \
\begin{pmatrix}
\phi^1_x(l_1) \\
\psi^1_x(l_1)
\end{pmatrix}=
\{M_{11}(\lambda)\}
\begin{pmatrix}
\zeta \\
\nu
\end{pmatrix}.
\end{eqnarray*}
In addition, the
 Cauchy problems on the edges $e_2,\ldots,e_{n-1}$ are solved:
\begin{eqnarray*}
-\phi_{xx}^i=\lambda \varrho_i(x)\phi^i,\quad
-\psi_{xx}^i=\lambda \mu_i(x)\psi^i,\quad x\in (0,l_i) \\
\begin{pmatrix}
\phi^1(l_i) \\
\psi^1(l_i)
\end{pmatrix}=
\begin{pmatrix}
0 \\
0
\end{pmatrix}, \ \
\begin{pmatrix}
\phi^i_x(l_i) \\
\psi^i_x(l_i)
\end{pmatrix}=
\{M_{1i}(\lambda)\}
\begin{pmatrix}
\zeta \\
\nu
\end{pmatrix}.
\end{eqnarray*}
Since we have already found the parameters of the equations on the
edges $e_i,$ the matrices $D_i$ and the angles $\alpha_{in}$ for
$i=1,\ldots, n-1$, the function $\{\Phi,\Psi\}$  is known on these
edges. At the internal vertex the compatibility conditions read
\begin{equation}
\begin{pmatrix}
\phi^1(0,\lambda)\\
\psi^1(0,\lambda)
\end{pmatrix}=S_{1n}\begin{pmatrix}
\phi^n(0,\lambda)\\
\psi^n(0,\lambda)
\end{pmatrix}, \ \
\sum_{j=1}^{n-1} S_{1j}D_j\begin{pmatrix}
\phi^j_x(0,\lambda)\\
\psi^j_x(0,\lambda)
\end{pmatrix}+S_{1n}D_n\begin{pmatrix}
\phi^n_x(0,\lambda)\\
\psi^n_x(0,\lambda)
\end{pmatrix}=0. \label{ffi5}
\end{equation}
Using these conditions and the definition of  $\widetilde M_{00},$
\begin{equation*}
\begin{pmatrix}
\phi^n_x(0,\lambda)\\
\psi^n_x(0,\lambda)
\end{pmatrix}=\widetilde M_{00}(\lambda)\begin{pmatrix}
\phi^n(0,\lambda)\\
\psi^n(0,\lambda)
\end{pmatrix},
\end{equation*}
we get the equations
\begin{equation}
\label{m2} \sum_{j=1}^{n-1} S_{1j}D_j \begin{pmatrix}
\phi^j_x(0,\lambda)\\
\psi^j_x(0,\lambda)
\end{pmatrix}
+ S_{1n}D_n\widetilde M_{00}(\lambda)(S_{1n})^{-1}
\begin{pmatrix}
\phi^1(0,\lambda)\\
\psi^1(0,\lambda)
\end{pmatrix}=0.
\end{equation}
Choosing two  linearly independent  vectors  $\begin{pmatrix}
\zeta\\
\nu
\end{pmatrix}$
 in (\ref{Star_BCN}), we get the corresponding vectors $\begin{pmatrix}
\phi^1(0,\lambda)\\
\psi^1(0,\lambda)
\end{pmatrix}$ in (\ref{m2}) to be linearly independent.
 Therefore, equation  (\ref{m2}) determines $\widetilde M_{00}(\lambda)$.

 Using the connection of the dynamical and spectral
data (\ref{Sp_Dyn}), we can recover the $\widetilde R_{00},$ the
component of the response function associated with the new boundary
node of $\widetilde\Omega,$ and reduce our problem to the inverse
problem for one edge.

We combine all results of this section in the next result.
\begin{thm}
Let $\Omega$ be the a star graph consisting of $n$ edges. Then the
graph and the parameters of the systems (\ref{Star_sp_eqn}),
(\ref{Star_cont}) and (\ref{Star_Kirhg}) are determined by the
reduced TW matrix $\left\{M_{ij}(\lambda)\right\}$, $1\leqslant
i,j\leqslant n-1$.
\end{thm}
\begin{remark} \label{mr2}
 Taking into account the connection between the response operator and TW matrix
 function (as in equation \eqref{Sp_Dyn}) we can claim that our inverse problem is solvable if we know
 the reduced response matrix function  $\textbf{R}(t) = \{R_{ij}(t)\},$
 $1\leqslant
i,j\leqslant n-1,$ for all $t>0.$ Moreover, we can recover the tree
and all parameters of the system from   $\textbf{R}(t), \, t \in
(0,T] ,\,$ for $T$ being finite but large enough. More precisely,
$T$ must be greater than or equal to the controllability time of the
system.  The sharp controllability time of system \eqref{Star_eqn}--
\eqref{Star_Kirhg} is equal to $2 L^*_n, \; L^*_n:=\max\{L_n,M_n\}
+\max_{i=1,\ldots,n-1} \max\{L_i,M_i\}$; see \cite{LLS1994}. We
discuss  connections between the dynamical and spectral versions of
the inverse problem in more detail for the general tree in the end
of Section 7.

If the whole response operator
$\left\{R_{ij}(t)\right\}$, $1\leqslant
i,j\leqslant n,$ known for $t \in [0,\, 2 L^*], \; L^*:=\max_{i=1,\ldots,n} \max\{L_i,M_i\},$
our inverse problem can clearly be solved without Step 5
 (and with the reduced Step 4).
\end{remark}

\section{Inverse problem on an arbitrary tree}
Let $\Omega$ be a finite tree with $m$ boundary   points
$\Gamma=\{v_1,\ldots,v_m\} $. Any boundary vertex of the tree can be
taken as a root; therefore  without  loss of generality, we can
assume that the boundary vertex $v_m$ is a root of the tree. We put
$\Gamma_m=\Gamma \setminus \{v_m\}$ and consider the dynamical
problem $\mathbf{D}$ and the spectral problem $\mathbf{S}$ on
$\Omega$. The reduced response function ${\bf
R}(t)=\{R_{ij}(t)\}_{i,j=1}^{m-1}$ and the TW matrix ${\bf
M}(\lambda)=\{M_{ij}(\lambda)\}_{i,j=1}^{m-1}$ associated with
 boundary points from $\Gamma_m$ are constructed as in
Section 3 and serve as data for our inverse problem.

{\bf Step 1.} Identifying the edges connected at the same vertex.

Using the diagonal elements  $R_{ii}(t)$,
$i=1,\ldots,m-1$ we can recover the densities $\varrho_i,$ $\mu_i$ on
the corresponding boundary edges and the optical lengths of the channels $L_i,$ $M_i$.
Two boundary edges, $e_i$ and $e_j,$  have a common
vertex if and only if
\begin{equation}
\{R_{ij}\}_{11}(t)=\left\{\begin{array}l =0,\quad \text{for}\,\, t<L_i+L_j\\
\not=0,\quad \text{for}\,\, t>L_i+L_j
\end{array}\right.,\quad 1\leqslant i,j\leqslant m-1.
\end{equation}
This relation allows us to divide the boundary edges into groups,
such that edges from one group have a common vertex. We call these
groups pre-sheaves.  More exactly, we introduce the following

{\bf Definition.}   We consider a subgraph of $\Omega$ which is a
star graph
consisting of {\sl all} edges incident to an internal vertex $v.$
 This star graph is called a \emph{pre-sheaf}
if it contains at least one  boundary edge of  $\Omega.$ A pre-sheaf is called a
\emph{sheaf}
if all but one its edges are the boundary edges of  $\Omega.$

The sheaves are especially important to our identification
algorithm. To extract them we denote the found pre-sheaves by $P_1,
\ldots, P_N$, and define  the optical distance $d(P_k,P_n)$ between
two pre-sheaves in the following way. We take boundary edges $e_i\in
P_k$ and $e_j\in P_n$ with the optical lengths of their first
channels  $L_i$ and $L_j$; we then  put
\begin{equation*}
d(P_k,P_n)=\max\{t>0\, : \,
\{R_{ij}\}_{11}(t)-L_i-L_j=0\}.
\end{equation*}
Clearly this definition does not depend on the particular choice
of $e_i,e_j$ and gives the distance between the internal vertices of the
 pre-sheaves
$P_k $ and $ P_n.$ Then we consider
\begin{equation*}
\max_{k,n\in 1,\ldots,N,\,k\not=n}d(P_k,P_n).
\end{equation*}
It is not difficult to see that two pre-sheaves on which this
maximum is attained (denoted here as $P$ and $P'$) are sheaves.
Indeed, since $\Omega$ is a tree, there is only one path between $P$
and $P'$. If we assume the existence of  an ``extra" internal edge
in $P$ or $P'$, this leads to contradiction, since there would
necessarily exist sheaves with a distance between them  greater than
$d(P,P')$.

 Having extracted a sheaf, we proceed with our identification procedure.

We consider now a sheaf consisting, say, of the boundary vertices
$v_1,\ldots,v_{m_0}$ from $\Gamma_m,$ the corresponding boundary
edges $e_1,\ldots,e_{m_0}$ and an internal edge $e_{m'_0}.$   We
identify each edge $e_i, \, i=1,\ldots,m_0$, with the interval
$[0,l_i]$ and the vertex $v_{m'_0},$  the internal vertex of the
sheaf, ---  with the set of common endpoints $x=0.$ At this point it
is convenient to renumerate the edge $e_{m'_0}$ as $e_{0}$ and the
vertex  $v_{m'_0}$ as  $v_{0}$. Applying the techniques of Section 4
(Steps 1--4) we find the  densities and lengths $\varrho_i,$
$\mu_i$, $l_i$, $i=1,\ldots,m_0$, angles $\alpha_{ij}$,
$i,j=0,1,\ldots,m_0,$ and matrices $D_i,$ $i=0,1,\ldots,m_0.$

{\bf Step 2.} The leaf-peeling method.

We consider a new tree, $\widetilde\Omega,$  which is $\Omega$
without the edges $ e_i$ and vertices $v_i, \, i=1,\ldots,m_0$. Our
goal in this step is to find $\widetilde {\bf M}(\lambda),$ the
(reduced) TW matrix associated with $\widetilde\Omega$. As in the
previous sections we identify the edges $e_i,$  $i=1,\ldots,m_0$,
with the intervals $(0,l_i)$ and the boundary vertices  $v_i$ ---
with  the interval endpoints $l_i.$

 By $\{\Phi,\Psi\}$ we
denote the solution to (\ref{2.7}), (\ref{2.8}) satisfying the
following boundary conditions
\begin{equation}
\label{TwoStGen} \{\Phi,\Psi\}=\{\zeta,\nu\},\,\,\text{at
$v_1$},\quad \{\Phi,\Psi\}=\{0,0\}, \,\,\text{at $v_i,\,\,
2\leqslant i \leqslant m$}.
\end{equation}
Note that the solution to (\ref{2.7}), (\ref{2.8}), (\ref{TwoStGen})
on the edge $e_1$
 solves the following Cauchy problem  on the edge $e_1:$
\begin{eqnarray*}
-\phi_{xx}^1=\lambda \varrho_1(x)\phi^1,\quad
-\psi_{xx}^1=\lambda \mu_1(x)\psi^1,\quad x\in (0,l_1), \\
\begin{pmatrix}
\phi^1(l_1) \\
\psi^1(l_1)
\end{pmatrix}=
\begin{pmatrix}
\zeta \\
\nu
\end{pmatrix}, \ \
\begin{pmatrix}
\phi^1_x(l_1) \\
\psi^1_x(l_1)
\end{pmatrix}=
\{M_{11}(\lambda)\}
\begin{pmatrix}
\zeta \\
\nu
\end{pmatrix}
\end{eqnarray*}
and also the Cauchy problems on the edges $e_2,\ldots,e_{m_0}:$
\begin{eqnarray*}
-\phi_{xx}^i=\lambda \varrho_i(x)\phi^i,\quad
-\psi_{xx}^i=\lambda \mu_i(x)\psi^i,\quad x\in (0,l_i) \\
\begin{pmatrix}
\phi^1(l_i) \\
\psi^1(l_i)
\end{pmatrix}=
\begin{pmatrix}
0 \\
0
\end{pmatrix}, \ \
\begin{pmatrix}
\phi^i_x(l_i) \\
\psi^i_x(l_i)
\end{pmatrix}=
\{M_{1i}(\lambda)\}
\begin{pmatrix}
\zeta \\
\nu
\end{pmatrix}.
\end{eqnarray*}
Thus, the function $\{\Phi,\Psi\}$ and its derivative are known on
the edges $e_1,\ldots,e_{m_0}$. At the internal vertex $v_{0}$
the KN conditions hold:
\begin{eqnarray}
\begin{pmatrix}
\phi^1(0,\lambda)\\
\psi^1(0,\lambda)
\end{pmatrix}=S_{10}\begin{pmatrix}
\phi^{0}(0,\lambda)\\
\psi^{0}(0,\lambda)
\end{pmatrix}, \ \
\sum_{j=1}^{m_0}S_{1j}D_j\begin{pmatrix}
\phi^j_x(0,\lambda)\\
\psi^j_x(0,\lambda)
\end{pmatrix}
+ S_{10}D_{0}\begin{pmatrix}
\phi^{0}_x(0,\lambda)\\
\psi^{0}_x(0,\lambda)
\end{pmatrix}=0. \label{Con1}
\end{eqnarray}
Using these conditions and the definition of the component of the
(reduced) TW matrix of the new graph $\widetilde \Omega$ associated
with the edge $e_0$
\begin{equation*}
\begin{pmatrix}
\phi^{0}_x(0,\lambda)\\
\psi^{0}_x(0,\lambda)
\end{pmatrix}=\widetilde M_{{0}{0}}(\lambda)\begin{pmatrix}
\phi^{0}(0,\lambda)\\
\psi^{0}(0,\lambda)
\end{pmatrix},
\end{equation*}
we obtain
\begin{eqnarray}  \label{m3}
\sum_{j=1}^{m_0} S_{1j}D_j \begin{pmatrix}
\phi^j_x(0,\lambda)\\
\psi^j_x(0,\lambda)
\end{pmatrix}
+S_{{1{0}}}D_{0}\widetilde
M_{{0}{0}}(\lambda)(S_{10})^{-1}
\begin{pmatrix}
\phi^1(0,\lambda)\\
\psi^1(0,\lambda)
\end{pmatrix}=0.
\end{eqnarray}
Equation (\ref{m3}) determines the matrix $\widetilde
M_{{0}{0}}(\lambda)$. By the definition of the reduced TW matrix
$\bf{M}(\lambda),$  we have
\begin{equation*}
\{M_{1j}(\lambda)\}\begin{pmatrix} \zeta\\
\nu
\end{pmatrix}=
\begin{pmatrix}
\phi^j_x(v_j)\\
\psi^j_x(v_j)
\end{pmatrix},\quad  j=m_0+1, \ldots, m-1.
\end{equation*}
On the other hand, by the definition of the reduced TW matrix for
the new tree $\widetilde\Omega,$
\begin{equation*}
\begin{pmatrix}
\phi^j_x(v_j)\\
\psi^j_x(v_j)
\end{pmatrix}=
\{\widetilde M_{0j}(\lambda)\}\begin{pmatrix} \phi^{0}(0)\\
\psi^{0}(0)
\end{pmatrix},\quad j=m_0+1, \ldots, m-1.
\end{equation*}
Thus, the component $\{\widetilde M_{0j}(\lambda)\}$ of the TW
matrix can be found from the equation; we use also the first
equality of (\ref{Con1})
\begin{equation} \label{ml1}
\{\widetilde M_{0j}(\lambda)\} \left(S_{10}\right)^{-1}
\begin{pmatrix}
\phi^1(0)\\
\psi^1(0)
\end{pmatrix}=
\{M_{1j}(\lambda)\}\begin{pmatrix} \zeta\\
\nu
\end{pmatrix},\quad  j=m_0+1, \ldots, m-1.
\end{equation}

In order to find the components $\widetilde M_{i0}(\lambda)$,
$i=m_0+1, \ldots, m-1$, we fix $i$ and denote by $\{\Phi,\Psi\}$ the solution
to (\ref{2.7}), (\ref{2.8}) with the boundary conditions
\begin{equation}
\label{TwoStGen1} \{\Phi,\Psi\}=\{\zeta,\nu\},\,\,\text{at
$v_i$},\quad \{\Phi,\Psi\}=\{0,0\}, \,\,\text{at $v_j,\,\,
j=1,\ldots, m,\,j\not= i $}.
\end{equation}
Note that on the edges $e_j,\, j=1, \ldots, m_0,$  this function
satisfies the equations
\begin{eqnarray*}
-\phi_{xx}^j=\lambda \varrho_j(x)\phi^j,\quad
-\psi_{xx}^j=\lambda \mu_j(x)\psi^j,\quad x\in (0,l_j), \\
\begin{pmatrix}
\phi^j(l_j) \\
\psi^j(l_j)
\end{pmatrix}=
\begin{pmatrix}
0 \\
0
\end{pmatrix}, \ \
\begin{pmatrix}
\phi^j_x(l_j) \\
\psi^j_x(l_j)
\end{pmatrix}=
\{M_{ij}(\lambda)\}
\begin{pmatrix}
\zeta \\
\nu
\end{pmatrix}.
\end{eqnarray*}
Thus, the function $\{\Phi,\Psi\}$ and its derivative are known on
the edges $e_1,\ldots,e_{m_0}$. Using conditions
(\ref{Con1})we can find the vectors $\begin{pmatrix}
\phi^{0}(0,\lambda)\\
\psi^{0}(0,\lambda)
\end{pmatrix}
$, $\begin{pmatrix}
\phi^{0}_x(0,\lambda)\\
\psi^{0}_x(0,\lambda)
\end{pmatrix}$.
We emphasize that the function $\{\Phi,\Psi\}$ does not satisfy zero
Dirichlet conditions at $v_{0}.$ We can obtain, however,  a function
satisfying this condition if we subtract from $\{\Phi,\Psi\}$ a
solution  $\{\tilde\Phi,\tilde\Psi\}$  of (\ref{2.7}), (\ref{2.8})
on $\tilde{\Omega}$ with the boundary conditions
\begin{equation*}
 \{\tilde\Phi,\tilde\Psi\}=\{\phi^{0}(0),\psi^{0}(0\},\,\,\text{at
$v_{0}$},\quad \{\tilde\Phi,\tilde\Psi\}=\{0,0\}, \,\,\text{at $v_j,\,\,
j=m_0+1,\ldots, m $}.
\end{equation*}
Therefore, the entries $\widetilde M_{i0}(\lambda)$,
$i=m_0+1,\ldots,m-1$ can be obtained from the equations
\begin{equation} \label{ml2}
\begin{pmatrix}
\phi^{0}_x(0,\lambda)\\
\psi^{0}_x(0,\lambda)
\end{pmatrix}-\widetilde M_{00}(\lambda)\begin{pmatrix}
\phi^{0}(0,\lambda)\\
\psi^{0}(0,\lambda)
\end{pmatrix}=\widetilde M_{i{0}}(\lambda)\begin{pmatrix}
\zeta\\
\nu
\end{pmatrix}.
\end{equation}
Correspondingly, the  entries $\widetilde M_{ij}(\lambda)$,
$i,j=m_0+1,\ldots,m-1$ can be found from the equations
\begin{equation}  \label{ml3}
M_{ij}(\lambda)\begin{pmatrix}
\zeta\\
\nu
\end{pmatrix}-\widetilde M_{0j}(\lambda)\begin{pmatrix}
\phi^{0}(0,\lambda)\\
\psi^{0}(0,\lambda)
\end{pmatrix}=\widetilde M_{ij}(\lambda)\begin{pmatrix}
\zeta\\
\nu
\end{pmatrix}.
\end{equation}
The  described procedure reduces the original inverse problem to the
inverse problem on a smaller subgraph.  This procedure may be
continued and, since the graph $\Omega$ is finite, it ends after a
finite number of steps. We proved the next result.
\begin{thm} \label{trs}
Let $\Omega$ be an arbitrary tree. Then the tree and the parameters
of system  (\ref{2.7}), (\ref{2.8}) are determined by the  entries
($2\times 2$ matrices) $M_{ij}(\lambda),$ $1\leqslant i,j\leqslant
m-1$ of the reduced TW matrix function.
\end{thm}

To formulate the result concerning the dynamical inverse problem, we
define the optical distance between the vertex $v_m$ and the rest
part of the boundary, $\Gamma_m:$
$$ \tau(v_m,\Gamma_m)=\max_{i=1, \ldots, m-1} \sum_{j:\, e_j \in \pi_i} \max \{L_j,M_j\}, $$
where $\pi_i$ is the path connecting $v_i$ with $v_n.$
\begin{thm}  \label{trd}
The tree $\Omega$ and  the
parameters of system  (\ref{2.10}) are
determined by the  entries  ($2\times 2$ matrices) $R_{ij}(t),$
$1\leqslant i,j\leqslant m-1,$ of the reduced response matrix function ${\bf R}(t)$ known for $t \in [0,2 \tau(v_m,\Gamma_m)].$
\end{thm}
We give a sketch of the proof; a detailed investigation of the
relation between spectral and dynamical inverse problems for
two-velocity systems on graphs will be presented in a separate
paper. Exact controllability of the system (\ref{2.10}) in the time
interval  $ [0,2 \tau(v_m,\Gamma_m)]$ was proved in \cite{LLS1994}.
One of the important results of the BC method allows us to recover
the spectral data of a system from the response operator, provided
the system is controllable and  the response operator is known on
the corresponding time interval \cite{ABR,ALP,B96}. Here under the
spectral data we mean the eigenvalues of the problem
(\ref{2.7})--(\ref{2.9}) and the derivatives
 of the eigenfunctions on the boundary $\Gamma.$ It is a fundamental fact of the spectral theory of ordinary differential operators
that the spectral data uniquely determines the TW function and vice
versa;
 see, e.g. \cite{AMR}. This result can be extended to
equations on graphs \cite{AB}. More specifically, we can claim that
the reduced response matrix function ${\bf R}(t)$ known on the
controllability time interval of a system determines the reduced TW
matrix function  (see, e.g. \cite{AB,AB1}). Therefore, Theorem
\ref{trd} follows from Theorem \ref{trs}.

\begin{remark}
Our method (as a version of the BC and leaf-peeling methods)
 leads to a stable algorithm that can be implemented numerically.
  For the (one-velocity) wave equation on a star graph,
numerical experiments were presented in \cite{ABM} with the
dynamical response operator as inverse data. We plan now to organize
numerical experiments for the two-velocity wave equation on graphs.
\end{remark}

\begin{acknowledgement}
Sergei Avdonin and Victor Mikhaylov wish to thank the Instituto de
F\'isica y Matem\'aticas, Universidad Michoacana de San Nicol\'as de
Hidalgo in Morelia, M\'exico for warm hospitality during their short
visits to this institute, and  acknowledge the support of
 CONACyt--M\'exico under the grant No.153184. Abdon Choque Rivero was
 also supported by the same CONACyt grant.
Sergei Avdonin  was partly supported by the NSF grant DMS 1411564,
Guenter Leugering was partly supported by the DFG-EC 315
  ''Engineering of advanced materials", and Victor Mikhaylov was supported by
following grants: NSh-1771.2014.1, RFBR 14-01-31388, and NIR SPbGU
11.38.263.2014.

The authors are very thankful to anonymous referees whose remarks
helped to improve the presentation of the material.
\end{acknowledgement}


\begin{thebibliography}{[1]}
\bibitem{A}
  S. Avdonin,  Control problems on quantum
graphs, in: {Analysis on Graphs and its Applications}, P. Exner, J.
Keating, P. Kuchment, T. Sunada, A. Teplyaev (eds.), Proceedings of
Symposia in Pure Mathematics, AMS {\bf 77}, 507--521 (2008).

\bibitem{ABM}
{S. Avdonin, B. Belinskiy and M. Matthews}, {Dynamical inverse
problems on a metric tree}, Inverse Problems  {\bf 27}(7), 1--21
(2011).


\bibitem{ABR} {S.A. Avdonin, M.I. Belishev, and
Yu.S. Rozhkov}, { The BC--method in the inverse problem for the heat
equation}, J. Inverse and Ill-Posed Problems {\bf 5}, 309--322
(1997).

\bibitem{AB} {S. Avdonin and
J. Bell}, { Determining a distributed parameter in a neural cable
model via a boundary control method,} J. Mathematical Biology {\bf
67}(1), 123--141 (2013).

\bibitem{AB1} {S. Avdonin and
J. Bell}, {Determining a Distributed Conductance Parameter for a
Neuronal Cable Model Defined on a Tree Graph,}
Journal of Inverse Problems and Imaging,  submitted.

\bibitem{AvdoninIvanov1995} S.A.~Avdonin and S.A.~Ivanov,   Families
of Exponentials. The Method of Moments in Controllability Problems
for Distributed Parameter Systems (Cambridge Univ. Press., 1995).

\bibitem{AK}  S. A. Avdonin and P. B. Kurasov,  Inverse problems for quantum trees,
Inverse Probl. Imag. \textbf{2}(1), 1--21 (2008).

\bibitem{ALP} {S. Avdonin, S. Lenhart, and V.
Protopopescu}, {Solving the dynamical inverse problem for the
Schr\"{o}dinger equation by the Boundary Control method}, Inverse
Problems {\bf 18}, 41--57  (2002).

\bibitem{ALM} S. Avdonin, G. Leugering and V. Mikhaylov,  On an inverse
problem for tree-like networks of elastic strings, Zeit. Angew.
Math. Mech. {\bf 90}, 136--150 (2010).



\bibitem{AM} S. Avdonin and V. Mikhaylov,  Controllability of partial
 differential equations on
graphs, Appl. Math. {\bf 35}, 379--393 (2008).

\bibitem{AMm} S. Avdonin and V. Mikhaylov,
The boundary control approach to inverse spectral theory, Inverse
Problems {\bf 26}(4), 1--19 (2010).


\bibitem{AMR} S.\,A.  Avdonin, V.\,S. Mikhaylov and A.\,V. Rybkin,
 The boundary control approach to the Titchmarsh--Weyl
$m-$function, Comm. Math. Phys. {\bf 275}(3), 791--803 (2007).

\bibitem{B96} Belishev M.I., Canonical model of a
dynamical system with boundary control in inverse problem for the
heat equation,  { St. Petersburg Math. Journal} {\bf 7}(6), 869--890
(1996).

\bibitem{BIP07}
M.\,I. Belishev, Recent progress in the boundary control method.
Inverse Problems {\bf 23}(5), R1--R67 (2007).

\bibitem{B} M.\,I. Belishev,  Boundary spectral inverse problem on a class of
graphs (trees) by the BC method, Inverse Problems {\bf 20}, 647--672
(2004).

\bibitem{B1}  M.\,I. Belishev, Boundary control and inverse problems:
 one-dimensional variant of BC-method, Zap. Nauchn. Sem. S.-Peterburg.
Otdel. Mat. Inst. Steklov. (POMI) {\bf 354} (2008), Mat. Vopr. Teor.
Rasprostr. Vol. 37, 19--80; translation in J. Math. Sci. (N. Y.)
{\bf 155}(3), 343--378 (2008).

\bibitem{BM} M.\,I. Belishev and V.\,S. Mikhaylov,  Unified approach to classical equations of
inverse problem theory, J. Inverse Ill-Posed Probl. {\bf 20}(4),
461--488 (2012).

\bibitem{BV} M.\,I. Belishev and A.\,F. Vakulenko,  Inverse problems
on graphs: recovering the tree of strings by the BC-method, Inverse
Ill-posed Probl. Ser. {\bf 14}, 29--46 (2006).


\bibitem{BW} B.\,M. Brown and R. Weikard, A Borg-Levinson
theorem for trees, Proc. R. Soc. Lond. Ser. A Math. Phys. Eng. Sci.,
\textbf{461}(2062),  3231--3243 (2005).


\bibitem{DagerZuazua2006}
R. D{\'a}ger and E.~Zuazua, Wave propagation, observation and
control in 1-$d$ flexible multi-structures. Series:
Math{\'e}matiques \& Applications 50. Springer-Verlag, Berlin, 2006.

\bibitem{FY}
G. Freiling and V. Yurko,
 Inverse problems for Strum-Liouville operators on noncompact trees,
 Results Math. \textbf{50}(3--4), 195--212 (2007).

\bibitem{KKLM} A. Kachalov, Y. Kurylev, M. Lassas and N. Mandache, 
Equivalence of time-domain inverse problems and boundary spectral
problems, Inverse Problems {\bf 20}, 491--436 (2004).

\bibitem{ks1}
{V.~Kostrykin and  R.~Schrader},  Kirchhoff's rule for quantum
wires, J. Phys A: Math. Gen. {\bf 32}, 595--630 (1999).

\bibitem{kuch}  P.~Kuchment,  Quantum graphs: an introduction and
a brief survey, in: { Analysis on Graphs and its Applications}, P.
Exner, J. Keating, P. Kuchment, T. Sunada, A. Teplyaev (eds.),
Proceedings of Symposia in Pure Mathematics, AMS {\bf 77}, 219--314
(2008).

\bibitem{KuSt}
P.~Kurasov and F.~Stenberg,  On the inverse scattering problem on
branching graphs, J. Phys. A: Math. Gen. {\bf 35},  101--121 (2002).

\bibitem{LagneseLeugering2004}
J.\,E. Lagnese and  G. Leugering,  Domain decomposition methods in
optimal control of partial differential equations, ISNM.
International Series of Numerical Mathematics {\bf 148}, Basel:
Birkh\"auser, 2004.

\bibitem{LagneseLeugeringSchmidt1993}
J.\,E. Lagnese, G. Leugering and E.\,J.\,P.\,G.~Schmidt,  Control of
planar networks of Timoshenko beams. SIAM J. Control Optimization
{\bf 31}(3), 780--811 (1993).

\bibitem{LagneseLeugeringSchmidt1994}
J.\,E.~Lagnese, G. Leugering and E.\,J.\,P.\,G. Schmidt, On the
analysis and control of hyperbolic systems associated with vibrating
networks. Proc. R. Soc. Edinb., Sect. A {\bf 124}(1), 77--104
(1994).

\bibitem{LLS1994}
J.\,E. Lagnese, G. Leugering and E.\,J.\,P.\,G  Schmidt, { Modeling,
analysis and control of dynamic elastic multi-link structures},
Birkh{\"a}user Boston, Systems and Control: Foundations and
Applications, 1994.

\bibitem{LeugeringSchmidt2011}
G. Leugering,  and E.\,J.\,P.\,G. Schmidt,
  On exact controllability of networks of nonlinear elastic
              strings in 3--dimensional space,
 {Chin. Ann. Math. Ser. B} {33}(1), 33--60 (2012).

\bibitem{LeugeringSokolowski2009}
G. Leugering and J. Sokolowski, Topological sensitivity analysis for
elliptic problem on graphs, Control and Cybernetics {\bf 37 }(4),
971--997 (2008).



\bibitem{LeugeringZuazua1999}
G.~Leugering and E.~Zuazua, On exact controllability of generic
trees, ESAIM, Proc. {\bf 8}, 95--105 (2000).

\bibitem{Y} {V. Yurko}, {Inverse Sturm-Liouville operator on
graphs}, Inverse Problems {\bf 21},  1075--1086 (2005).
\end{thebibliography}
\end{document}